\documentclass[a4paper,11pt]{article}
\usepackage[utf8x]{inputenc}
\usepackage{amsthm,amsmath,amssymb,oldgerm}
\usepackage[margin=3cm]{geometry}
\theoremstyle{plain}
\newtheorem{thm}{Theorem}

\newtheorem{cor}[thm]{Corollary}
\theoremstyle{definition}
 
\theoremstyle{definition}
\newtheorem{rem}[thm]{Remark}

\let\Im\relax
\DeclareMathOperator{\Im}{Im}
\DeclareMathOperator{\bk} {{\it{\mathcal{B}_{X}^{k}}}}
\DeclareMathOperator{\bkm} {{\it{\mathcal{B}_{M,\mathcal{L}}^{k}}}}
\DeclareMathOperator{\hkm} {{\it{K_{M,\mathcal{L}}^{k}}}}

\DeclareMathOperator{\hyp}{\mu_{hyp}}

\DeclareMathOperator{\shyp}{\mu_{shyp}}

\DeclareMathOperator{\vx}{\mathrm{vol_{\mathrm{hyp}}}}

\DeclareMathOperator{\lk}{{\it{\overline{\square}_{k}}}}
\let\ell\relax
\DeclareMathOperator{\ell}{\mathcal{L}}
\title{Heat kernel approach for sup-norm bounds for cusp forms of integral and half-integral weight}
{\small\author{Anilatmaja Aryasomayajula}}
\date{}
\begin{document}
\maketitle
\begin{abstract}
\noindent In this article, using the heat kernel approach from \cite{bouche}, we derive sup-norm bounds for cusp 
forms of integral and half-integral weight. Let $\Gamma\subset \mathrm{PSL}_{2}(\mathbb{R})$ be a cocompact 
Fuchsian subgroup of first kind. For $k\in\frac{1}{2}\mathbb{Z}$ (or $k\in 2\mathbb{Z}$), let 
$S^{k}(\Gamma)$ denote the complex vector space of weight-$k$ cusp forms. 
Let $\lbrace f_{1},\ldots,f_{j_{k}} \rbrace$ denote an orthonormal basis of $S^{k}(\Gamma)$. In this article, we show 
that as $k\rightarrow \infty,$ the sup-norm for  $\sum_{i=1}^{j_{k}}y^{k}|f_{i}(z)|^{2}$ is bounded by $O(k)$, 
where  the implied constant is independent on $\Gamma$.  Furthermore, using results from \cite{berman}, we extend these results to the case when $\Gamma$ is cofinite. 

\vspace{0.2cm}\noindent
\end{abstract}
\section{Introduction}\label{introduction}
The main aim of this article is to demonstrate the usefulness of the heat kernel approach to study 
automorphic forms. J.~Jorgenson and J.~Kramer have successfully applied heat kernel analysis to study 
automorphic forms, and derived optimal estimates for various arithmetic invariants. We use a slightly different 
approach, which is more geometric and derive similar estimates for integral weight cusp forms as in 
\cite{jk2}. We extend these estimates to half-integral weight cusp forms. 

\vspace{0.025cm}
\paragraph{Notation}
Let $\Gamma\subset \mathrm{PSL}_{2}(\mathbb{R})$ be a cocompact Fuchsian subgroup of first kind acting on the 
hyperbolic upper half-plane $\mathbb{H}$, via fractional linear transformations. Let $X$ denote the quotient 
space $\Gamma\backslash\mathbb{H}$, and let $g$ denote the genus of $X$. We assume that  $g>1$. 

\vspace{0.2cm}
We assume that $\Gamma$ has no elliptic elements, so the quotient space $X$ admits the 
structure of a Riemann surface.  

\vspace{0.2cm}
We denote the $(1,1)$-form corresponding to the hyperbolic metric of $X$, which is compatible with the 
complex structure on $X$ and has constant negative curvature equal to minus one, by $\hyp(z)$. Locally, for 
$z\in X$, it is given by
\begin{equation*}
 \hyp(z)= \frac{i}{2}\cdot\frac{dz\wedge d\overline{z}}{{\Im(z)}^{2}}.
\end{equation*} 
Let $\shyp(z)$ denote the rescaled hyperbolic metric $\hyp(z)\slash \vx(X)$, which measures the volume of $X$ to be 
one. 

\vspace{0.2cm}
For $k\in\frac{1}{2}\mathbb{Z}$ (or $2\mathbb{Z}$), let $S^{k}(\Gamma)$ denote the complex vector space of weight-$k$ cusp forms. 
Let $j_{k}$ denote the dimension of $S^{k}(\Gamma)$, and let $\lbrace f_{1},\ldots,f_{j_{k}} \rbrace$ denote 
an orthonormal basis of $S^{k}(\Gamma)$ with respect to the Petersson inner product.  Put 
\begin{align*}
\bk(z):=\sum_{i=1}^{j_{k}} y^{k}|f_{i}(z)|^{2}.
\end{align*}

With notation as above, for $k\in\frac{1}{2}\mathbb{Z}$ (or $2\mathbb{Z}$), we have the 
following estimate
\begin{align}\label{estimate1}
\lim_{k}\sup_{z\in X}\frac{1}{k}\bk(z)= O(1),
\end{align}
where the implied constant is independent of $\Gamma$. 

\vspace{0.2cm}
Let $\Gamma$ now be a cofinite Fuchsian subgroup without elliptic elements, which implies that 
$X=\Gamma\backslash\mathbb{H}$ admits the structure of a noncompact Riemann surface of finite volume. Then, 
for any $z\in X$ and $k\in\frac{1}{2}\mathbb{Z}$, we have 
\begin{align}\label{estimate2}
\limsup_{k}\frac{1}{k}\bk(z)\leq\frac{1}{8\pi};
\end{align} 
and for $k\in 2\mathbb{Z}$, we have 
\begin{align}\label{estimate3}
\limsup_{k}\frac{1}{k}\bk(z)\leq\frac{1}{2\pi}.
\end{align}  

\vspace{0.025cm}
\paragraph{Applications and existing results}
When $X$ is comoact, in proving estimate \eqref{estimate1}, for any $z\in X$ and $k\in\frac{1}{2}\mathbb{Z}$, we first arrive at the 
following estimate
\begin{align*}
\lim_{k}\frac{1}{k}\bk(z) =O(1),
\end{align*}
where the implied constant is independent of $\Gamma$. So, let 
\begin{align*}
\lim_{k}\frac{1}{k}\bk(z)=C.
\end{align*}
Hence, for $z\in X$ and $k\in\frac{1}{2}\mathbb{Z}$ (or $2\mathbb{Z}$) and $f$ any smooth function on $X$, 
we observe that
\begin{align*}
 \lim_{k}\frac{1}{Ck}\int_{X}f(z)\bk(z)\shyp(z)=\int_{X}f(z)\shyp(z),
\end{align*}
which proves a version of arithmetic quantum unique ergodicity conjecture in the average case for half-integral weight 
cusp forms. 

\vspace{0.2cm}
For $k\in\frac{1}{2}\mathbb{Z}$ and $N\in\mathbb{N}$, let $f$ be any weight-$k$ cusp form with respect to the arithmetic subgroup 
$\Gamma_{0}(4N)$. Furthermore, let $f$ be normalized with respect to the Petersson 
inner-product. Then, in \cite{kiral}, Kiral has derived the following estimate 
\begin{align*}
 \sup_{z\in Y_{0}(N)}y^{k}|f(z)|^{2}=O_{k,\varepsilon}\big(N^{\frac{1}{2}-\frac{1}{18}+\varepsilon}\big),
\end{align*}
for any $\varepsilon>0$. Using the above estimate, one can derive 
\begin{align*}
 \sup_{z\in Y_{0}(N)}\mathcal{B}_{Y_{0}(N)}(z)=O_{k,\varepsilon}\big(N^{\frac{1}{2}-\frac{1}{18}+\varepsilon}\big),
\end{align*}
for any $\varepsilon>0$.

\vspace{0.2cm}
For $k\in2\mathbb{Z}$, Friedman, Jorgenson, and Kramer have derived the estimates $O(k)$ and 
$O(k^{3\slash 2})$ for $\bk(z)$ in the cocompact and cofinite cases, respectively in \cite{jk2}. Extending their 
heat kernel techniques from \cite{jk1} to higher weights, they showed that the implied constant is independent of the 
Fuchsian subgroup $\Gamma$.  

\vspace{0.2cm}
For $k\in\frac{1}{2}\mathbb{Z}$, estimate \eqref{estimate1} shows the dependence of $\bk$ on $k$, and our 
result holds true for any arbitrary cocompact Fuchsian group. 

\vspace{0.2cm}
For $k\in 2\mathbb{Z}$, estimate \eqref{estimate1} gives an alternate proof for the estimate of Jorgenson and Kramer in the cocompact 
case. Our method is more geometric and relies on  a theorem of Bouche in \cite{bouche}, which is based on the 
micro-local analysis of the heat kernel. 

\vspace{0.2cm}
Estimate \eqref{estimate1} extends with notational changes to cusp forms of integral and half-integral weight 
associated to cocompact Hilbert modular groups and Siegel modular groups.  

\vspace{0.2cm}
It is possible to extend estimate \eqref{estimate1} to cofinite Fuchsian subgroups by studying Bouche's 
methods in a cuspidal neighborhood.
\section{Heat kernels on compact complex manifolds}
In this section, we recall the main results from \cite{bouche} and \cite{berman}, which we use in 
the next section. 

\vspace{0.2cm}
Let $(M,\omega)$ be a compact complex manifold of dimension $n$ with a Hermitian metric $\omega$. Let 
$\ell$ be a positive Hermitian holomorphic line bundle on $M$ with the Hermitian metric given by 
$\|s(z)\|^{2}_{\ell}:=e^{-\phi(z)}|s(z)|^{2}$, where $s\in\ell$ is any section, and $\phi(z)$ is 
a real-valued function defined on $M$. 

\vspace{0.2cm}
For any $k\in\mathbb{N}$, let $\lk:=(\overline{\partial}^{\ast}+\overline{\partial})^{2}$ denote the $\overline{\partial}$-Laplacian 
acting on smooth sections of the line bundle $\ell^{\otimes k}$. Let $\hkm(t;z,w)$ denote the smooth kernel of the 
operator $e^{-\frac{2t}{k}\lk}$. We refer the reader to p. 2 in \cite{bouche}, for the details regarding the 
properties which uniquely characterize the heat kernel $\hkm(t;z,w)$. When $z=w$, the heat kernel $\hkm(t;z,w)$ 
admits the following spectral expansion
\begin{align}\label{spectralexpn}
\hkm(t;z,w)=\sum_{n\geq 0} e^{-\frac{2t}{k}\lambda_{n}^{k}}\varphi_{n}(z)\otimes\varphi_{n}^{\ast}(w),
\end{align}
where $\lbrace\lambda_{n}^k\rbrace_{n\geq0}$ denotes the set of eigenvalues of $\lk$ 
(counted with multiplicities), and $\lbrace\varphi_{n}\rbrace_{n\geq0}$ denotes a set of associated 
orthonormal eigenfunctions. 

\vspace{0.2cm}
Let $\lbrace s_{i}\rbrace$ denote an orthonormal basis of $H^{0}(M,\mathcal{L}^{\otimes k})$. For any $z\in M$, 
the Bergman kernel is given by
\begin{align}\label{bkdefn}
\bkm(z):= \sum_{i}\| s_{i}(z)\|_{\mathcal{L}^{\otimes k}}^{2}.
\end{align}
For any $z\in M$ and $t\in\mathbb{R}_{>0}$, from the spectral expansion of the heat kernel $\hkm(t;z,w)$ described in equation \eqref{spectralexpn}, 
it is easy to see that
\begin{align}\label{hkbkreln}
\bkm(t;z)\leq \hkm(t;z,z)\quad \mathrm{and}\quad
\lim_{t} \hkm(t;z,z)=\bkm(t;z).
\end{align}
For $z\in M$, let $c_{1}(\ell)(z):=\frac{i}{2\pi}\partial\overline{\partial}\phi(z)$ denote the first Chern form 
of the line bundle $\ell$. Let $\alpha_{1},\ldots,\alpha_{n}$ denote the eigenvalues of $\partial\overline{\partial}\phi(z)$ at 
the point $z\in M$. Then, with notation as above, 
from Theorem 1.1 in \cite{bouche}, for any $z\in M$ and $t\in (0,k^{\varepsilon})$, for 
a given $\varepsilon>0$ not depending on $k$, we have
\begin{align}\label{boucheeqn1}
\lim_{k}\frac{1}{k^{n}}\hkm(t;z,z)=\prod_{j=1}^{n}\frac{\alpha_{j}}{(4\pi)^{n}\sinh(\alpha_{j}t)}, 
\end{align}
and the convergence of the above limit is uniform in $z$. 

\vspace{0.2cm}
Using equations \eqref{hkbkreln} and \eqref{boucheeqn1}, in Theorem 2.1 in \cite{bouche}, Bouche derived the 
following asymptotic estimate 
\begin{align}\label{boucheeqn2}
\lim_{k}\frac{1}{k^{n}}\bkm(z)= O\big(\mathrm{det}_{\omega}\big(c_{1}(\ell)(z)\big)\big),
\end{align}
and the convergence of the above limit is uniform in $z\in X$. 

\vspace{0.2cm}
When $M$ is a noncompact complex manifold, using micro-local analysis of the Bergman kernel, in \cite{berman}, 
Berman derived the following estimate
\begin{align}\label{bermaneqn}
\limsup_{k}\frac{1}{k^{n}}\bkm(z)\leq \mathrm{det}_{\omega}\big(c_{1}(\ell)(z)\big).
\end{align}
\section{Estimates of cusp forms}
Let notation be as in Section \ref{introduction}. We allow $\Gamma$ to be cofinite, i.e. $X$ can be a 
noncompact Riemann surface of finite volume. Let $\Omega_{X}$ denote the cotangent bundle over $X$. Then, for any $k\in2\mathbb{Z}$, cusp forms of 
weight $k$ with respect to $\Gamma$ are global section of the line bundle $\Omega_{X}^{\otimes k\slash 2}$. 
Furthermore, recall that for any $f\in\Omega_{X}$, i.e., $f$ a weight-$2$ cusp form, the Petersson metric on the 
line bundle $\Omega_{X}$ is given by 
\begin{align}\label{peterssonip1}
 \|f(z)\|_{\Omega_{X}}^{2}:=y^{2}|f(z)|^{2}. 
\end{align}
Let $\omega_{X}$ denote the line bundle of cusp forms of weight $\frac{1}{2}$ over $X$. Then, for any 
$k\in\frac{1}{2}\mathbb{Z}$, cusp forms of weight-$k$ with respect to $\Gamma$ are global section of the line 
bundle $\omega_{X}^{\otimes 2k}$. Furthermore, recall that for any 
$f\in\omega_{X}$, i.e., $f$ a weight-$\frac{1}{2}$ cusp form, the Petersson metric on the line 
bundle $\omega_{X}$is given by 
\begin{align}\label{peterssonip2}
 \|f(z)\|_{\omega_{X}}^{2}:=y^{1\slash 2}|f(z)|^{2}. 
\end{align}

\vspace{0.2cm}
\begin{rem}
For any $z\in X$ and $k\in 2\mathbb{Z}$, from the definition of the Bergman kernel $\mathcal{B}_{X,\Omega_{X}}^{k
\slash 2}(z)$ for the line bundle $\Omega_{X}^{\otimes k\slash 2}$ from equation \eqref{bkdefn}, we have
\begin{align*}
\mathcal{B}_{X,\Omega_{X}}^{k\slash 2}(z)=\bk(z). 
\end{align*}
Similarly, for any $z\in X$ and $k\in \frac{1}{2}\mathbb{Z}$, from the definition of the Bergman kernel 
$\mathcal{B}_{X,\Omega_{X}}^{2k}(z)$ for the line bundle $\omega_{X}^{\otimes 2k}$ from equation \eqref{bkdefn}, we have
\begin{align}\label{remeqn2}
\mathcal{B}_{X,\omega_{X}}^{2k}(z)=\bk(z). 
\end{align}
\end{rem}

\vspace{0.2cm}
\begin{thm}\label{thm1}
Let notation be as in Section \ref{introduction}, and we assume that $\Gamma$ is cocompact, i.e. $X$ is compact. 
Then, for $k\in\frac{1}{2}\mathbb{Z}$ (or $2\mathbb{Z}$), we have the following estimate
\begin{align*}
\lim_{k}\sup_{z\in X}\frac{1}{k}\bk(z)=O(1), 
\end{align*}
where the implied constant is independent of $\Gamma$.   
\begin{proof}
We prove the theorem for $k\in \frac{1}{2}\mathbb{Z}$, and the case for $k\in2\mathbb{Z}$ follows 
automatically with notational changes. For any $z\in X$,  observe that
\begin{align*}
 c_{1}(\omega_{X})(z)=-\frac{i}{2\pi}\partial\overline{\partial}\log\big(y^{1\slash 2}|f(z)|^{2}\big),
\end{align*}
where $f$ is any cusp form of weight $\frac{1}{2}$. From the above equation, we derive 
\begin{align*}
c_{1}(\omega_{X})(z) =\frac{i}{16\pi}\cdot\frac{dz\wedge\overline{z}}{ y^{2}}=\frac{1}{8\pi}\hyp(z),
\end{align*}
which shows that the line bundle $\omega_{X}$ is positive, and $\mathrm{det}_{\hyp}\big(c_{1}(\omega_X)(z)
\big)=\frac{1}{8\pi}$. So applying estimate \eqref{boucheeqn2} to the complex manifold $X$ with its natural Hermitian 
metric $\hyp$ and the line bundle $\omega_{k}^{\otimes k}$, and using equation \eqref{remeqn2}, we find
\begin{align*}
\lim_{k}\frac{1}{k}\bk(z)=\lim_{k}\frac{1}{k}\mathcal{B}_{X,\omega_X}^{2k}(z)=O
\bigg(2\mathrm{det}_{\hyp}\big(c_{1}(\omega_X)(z)\big)\bigg)=O(1).
\end{align*}
As the above limit convergences uniformly in $z\in X$, and as $X$ is compact, we have
\begin{align*}
\sup_{z\in X} \lim_{k}\frac{1}{k}\bk(z)=\lim_{k}\sup_{z\in X}\frac{1}{k}\bk(z)=O(1),
\end{align*}
which completes the proof of the theorem.
\end{proof}
\end{thm}

\vspace{0.2cm}
\begin{cor}
Let notation be as in Section \ref{introduction}, and we assume that $\Gamma$ is cofinite, i.e. $X$ is a 
noncompact Riemann surface of finite volume. Then, for $z\in X$ and $k\in\frac{1}{2}\mathbb{Z}$, we have 
\begin{align*}
\limsup_{k}\frac{1}{k}\bk(z\leq\frac{1}{8\pi};
\end{align*}
and for  $k\in 2\mathbb{Z}$, we have 
\begin{align*}
\limsup_{k}\frac{1}{k}\bk(z\leq\frac{1}{2\pi}.
\end{align*}
\begin{proof}
The proof of the theorem follows from estimate \eqref{bermaneqn}, and from similar arguments as in Theorem \ref{thm1}. 
 \end{proof}
\end{cor}
\paragraph{Acknowledgements}
The author would like to thank J.~Kramer and J.~Jorgenson for introducing him to the area of automorphic forms 
and heat kernels. The author would like to express his gratitude to T. Bouche for providing him with helpful 
references, and to Archana S. Morye, for many helpful discussions and remarks.
{\small{}}

\vspace{0.3cm}
{\small{Anilatmaja Aryasomayajula\\ {\it{anilatmaja@gmail.com}}  \\ 
Department of Mathematics, \hfill\\University of Hyderabad, \\Prof. C.~R.~Rao Road, Gachibowli,\\
Hyderabad, 500046, India}}
\end{document}